\documentclass{article}
\usepackage{amssymb}
\usepackage{amsmath}
\usepackage{epsf}
\usepackage{color}
\def \C {{\Bbb C\,}}

\def\N {{\Bbb N}}
\def \Z {{\Bbb Z}}

\def\boB{{\cal B}}

\def\boE{{\cal E}}

\def\boH{{\cal H}}

\def\boM{{\cal M}}
\def\cc{c_1}
\def\ccc{c_2}
\def\hh{{\bf h}}

\def\pp{{\bf p}}
\def\qq{{\bf q}}
\def\tt{{\bf t}}

\def\cqfd{\hfill$\Box$}

\def\disp{\displaystyle}

\newtheorem{theorem}{Theorem}
\newtheorem{proposition}{Proposition}
\newtheorem{definition}{Definition}
\newtheorem{lemma}{Lemma}
\newtheorem{corollary}{Corollary}

\newtheorem{remark}{Remark}
\newtheorem{claim}{Claim}
\newtheorem{hypothesis}{Hypothesis}
\newtheorem{example}{Example}
\newtheorem{convention}{Convention}
\begin{document}
\title{Opening infinitely many nodes}
\author{Martin Traizet}
\maketitle
{\em Abstract : we develop a theory of holomorphic differentials
on a certain class of non-compact Riemann surfaces obtained by
opening infinitely many nodes.}

\medskip

\section{Introduction}
\label{section-introduction}
This paper is about holomorphic 1-forms,
on a certain class of non compact Riemann surfaces.
The compact case is a classical subject. The space
of holomorphic 1-forms (also called holomorphic differentials)
on a compact Riemann surface $\Sigma$ has complex dimension
equal to the genus $g$ of $\Sigma$. Moreover, a holomorphic differential
$\omega$ can be defined by prescribing its integrals 
(also called periods) on a suitable set of cycles, namely the cycles $A_1,\cdots,A_g$
of a canonical homology basis.

\medskip

The non-compact case is mostly unexplored. In this case, the space of
holomorphic differentials on $\Sigma$ is infinite dimensional. For any practical
purpose, it is required to put a norm on this space. So what we
are interested in are Banach spaces of holomorphic differentials on $\Sigma$.
The question which we would like to answer is :

\medskip

{\em Can we define a holomorphic
differential by prescribing its periods,
and in which Banach space does it live ?}

\medskip

\begin{figure}
\begin{center}
\epsfxsize=10cm
\epsffile{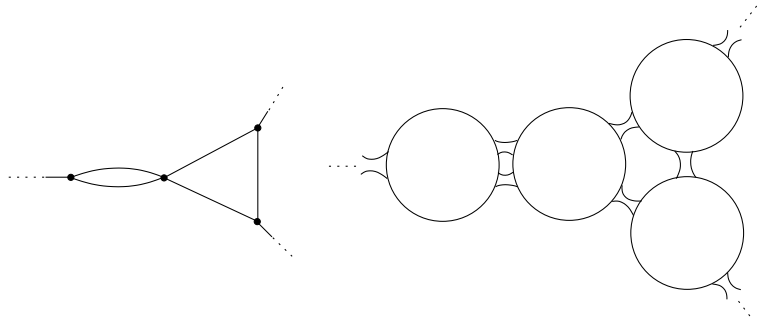}
\caption{A small portion of an infinite graph
and a picture of its associated Riemann surface.}
\label{fig1}
\end{center}
\end{figure}

We will answer this question on a class of non-compact Riemann surfaces which
are made of infinitely many Riemann spheres connected by small necks. More
precisely, consider an infinite graph $\Gamma$. Let $V$ denote its set of vertices
and $E$ its set of edges.
For each vertex $v\in V$, consider a Riemann sphere $S_v$. For any edge $e\in E$
from the vertex $v$ to the vertex $v'$, we connect $S_v$ and $S_{v'}$ by a small neck,
see figure \ref{fig1}.
We do this by opening nodes, a standard explicit construction.

\medskip

We denote by $\gamma_e$ the cycle around the neck corresponding to the edge $e$.
Essentially we prove that we can define a holomorphic differential $\omega$ by prescribing
its periods $\alpha_e$ on the cycles $\gamma_e$ for $e\in E$, provided the prescribed 
periods satisfy the obvious homological obstruction obtained from Cauchy theorem
in each Riemann sphere.

\medskip

One issue is to put norms on both the space of holomorphic differentials $\omega$ and the
space of period vectors $\alpha=(\alpha_e)_{e\in E}$. Our choice was to leave the later norm as general as possible,
and see what is the right norm on the space of holomorphic differentials so that
the operator which maps a holomorphic differential $\omega$ to its period vector $\alpha$ is
a Banach isomorphism. This result is theorem \ref{theorem2}.

\medskip

We also prove that the holomorphic differential $\omega$
obtained by prescribing its periods depends
smoothly (actually, analytically)
on all parameters which enter the construction of $\Sigma$. Roughly speaking, these parameters can be thought of as the position and size of the necks.

\medskip

Because our Riemann surface is not compact, one can ask what is the
asymptotic behavior of the holomorphic differential $\omega$.
This can be done using weighted norms and will be a recurrent question
in this paper.

\medskip

Finally, in the last section of the paper,
we consider meromorphic differentials. In the compact case, it is well known
that one can define a meromorphic differential by prescribing its poles, principal part at each
pole and periods. We prove such a result in our non compact setting
(a Mittag Leffler type result).

\medskip

The theory of holomorphic differentials on non compact Riemann surfaces which we develop in this
paper, besides being interesting in its own right, has applications to minimal surfaces theory,
via the classical Weierstrass Representation. It can be used to glue infinitely many minimal surfaces
together. See \cite{filippo} for an example.
\section{Opening nodes}
\label{section-openingnodes}
Consider a connected  graph $\Gamma$.
It may have multiple edges and edges from a vertex to itself (this
is usually called a multi-graph).
We denote by $V$ its set of vertices and $E$ its set of edges.
We are mostly interested in the case where $V$ is not finite.
We assume that $\Gamma$ is oriented, which means that each edge has an orientation.
The orientation is arbitrary and will only be used to orient certain curves.
If $e\in E$ is an oriented edge from vertex $v$ to vertex $v'$, we say that
$v$ is the starting point of $e$ and $v'$ is the endpoint.

\medskip

For a vertex $v\in V$, the set of edges having $v$ as starting point (resp. endpoint)
will be denoted $E_v^-$ (resp. $E_v^+$).
The set of edges adjacent to the vertex $v$ is denoted $E_v=E_v^-\cup E_v^+$.
The degree of $v$ is the cardinal of $E_v$, which we assume is always finite.
Note that together with the fact that $\Gamma$ is connected,
this implies that $V$ is countable.

\medskip

To each vertex $v\in V$ we associate a Riemann sphere $\C\cup\{\infty\}$ denoted $S_v$. 
For each oriented edge $e\in E$ from $v$ to $v'$,
we choose a point
$p_e^-$ in $S_v$ and a point $p_e^+$ in $S_{v'}$.
We assume that for each vertex $v$,
the points $p_e^-$ for $e\in E_v^-$ and $p_e^+$ for $e\in E_v^+$ are distinct.

\medskip

We consider the disjoint union of all $S_v$, $v\in V$, and
for each edge $e\in E$,
we identify the points $p_e^-$ and $p_e^+$.
This creates a node which we call $p_e$.
We call $\Sigma_0$ the resulting Riemann surface with nodes.

\medskip

We open nodes in the standard way as follows,
see for example \cite{imayoshi}.
We use the notation $D(a,R)$ for
the disk of center $a$ and radius $R$ in $\C$.
We consider, for each edge $e\in E$ from $v$ to $v'$,
some local complex coordinates
$z_e^-:D_e^-\subset S_v\stackrel{\sim}{\to} D(0,\rho)$ and
$z_e^+:D_e^+\subset S_{v'}\stackrel{\sim}{\to} D(0,\rho)$ in a neighborhood
of $p_e^-$ and $p_e^+$ respectively,
such that $z_e^-(p_e^-)=0$ and $z_e^+(p_e^+)=0$.
We assume that for each vertex $v$, the disks $D_e^-$ for $e\in E_v^-$
and $D_e^+$ for $e\in E_v^+$ are disjoint in $S_v$.

\medskip

Consider again the disjoint union of all Riemann spheres $S_v$, $v\in V$.
For each edge $e\in E$ from $v$ to $v'$, choose some complex number $t_e\in D(0,\rho^2)$.
If $t_e\neq 0$, remove the disk
$|z_e^-|\leq \frac{|t_e|}{\rho}$ from $D_e^-$ and
$|z_e^+|\leq \frac{|t_e|}{\rho}$ from $D_e^+$.
We identify each point $z\in D_e^-$ with the point $z'\in D_e^+$ 
such that
$z_e^-(z)z_e^+(z')=t_e.$
This creates a neck connecting $S_v$ and $S_{v'}$.
If $t_e=0$, then we identify $p_e^-$ and $p_e^+$ as before to create a node.

\medskip

Doing this for all edges defines a (possibly noded) Riemann surface which we call
$\Sigma_{\tt}$. Here ${\tt}$
denotes the sequence $(t_e)_{e\in E}$.
When all $t_e$ are nonzero, $\Sigma_{\tt}$ is a genuine Riemann surface.
It is compact when $\Gamma$ is a finite graph.
\section{Regular differentials}
\label{section-regular}
Regular differentials are the natural generalisation of holomorphic 1-forms to
Riemann surfaces with nodes.
\begin{definition}[Bers \cite{masur}]
\label{def-regulardifferential}
A differential $\omega$
on a Riemann surface with nodes $\Sigma$ is called {\em regular} if it
is holomorphic away from the nodes and for each node $p$ (obtained by identifying
$p^-$ and $p^+$), it has simples poles at $p^-$ and $p^+$ with opposite residues.
\end{definition}
For each edge $e\in E$, the boundary $\partial D_e^+$ is homologous 
in $\Sigma_{\tt}$ to
$-\partial D_e^-$. We define the cycle $\gamma_e$ as the homology class of
$\partial D_e^+$.
We have
\begin{equation}
\label{eq-homology}
\forall v\in V,\qquad
\sum_{e\in E_v^-} \gamma_e=\sum_{e\in E_v^+} \gamma_e.
\end{equation}
We want to define a regular differential $\omega$ on $\Sigma_{\tt}$ by prescribing
its periods on the cycles $\gamma_e$ :
$$\forall e\in E, \qquad \int_{\gamma_e} \omega = \alpha_e, \qquad \alpha_e\in \C$$
By Cauchy theorem and equation \eqref{eq-homology}, a necessary condition is
\begin{equation}
\label{eq-compatible}
\forall v\in V,\qquad
\sum_{e\in E_v^-} \alpha_e=\sum_{e\in E_v^+} \alpha_e.
\end{equation}
If $\Gamma$ is a finite graph, equation \eqref{eq-compatible} is the only
obstruction to define a regular differential on $\Sigma_{\tt}$ :
\begin{theorem}[Fay \cite{fay}]
\label{theorem-fay}
If $\Gamma$ is a finite graph, then for $\tt$ small enough, the operator 
$\omega\mapsto (\int_{\gamma_e}\omega)_{e\in E}$ is an isomorphism from
the space of regular differentials
on $\Sigma_{\tt}$ to the space of vectors $(\alpha_e)_{e\in E}$ which
satisfy the compatibility condition \eqref{eq-compatible}
\end{theorem}
\begin{remark}
Fay's theorem is more general, as he does not require the parts $S_v$ to
be spheres.
\end{remark}
\subsection{Admissible coordinates}
\label{section-admissiblecoordinates}
In the case of an infinite graph, we need to make some
assumptions on the coordinates used to open nodes.
We denote by $z$ the standard coordinate in each sphere $S_v$,
and
write $z_e^{\pm}$ to designate either $z_e^+$
or $z_e^-$.
For each $e\in E$, the ratio $\left|\frac{z_e^{\pm}}{z-p_e^{\pm}}\right|$ extends continuously at $p_e^{\pm}$ so
is bounded from above and below in $D_e^{\pm}$ by positive
numbers. We require these numbers to be independant of $e\in E$.
The following definition summarizes our requirements on the
coordinates.
\begin{definition}
\label{def-admissiblecoordinates}
We say that the coordinates $(z_e^{\pm})_{e\in E}$ are admissible
if
\begin{enumerate}
\item all points $p_{e}^-$ and $p_e^+$ are different from $\infty$,
\item for each $e\in E$,
$z_e^{\pm}:D_e^{\pm}\stackrel{\sim}{\to} D(0,\rho)$
is a diffeomorphism such that $z_e^{\pm}(p_e^{\pm})=0$
($\rho$ is independant of $e$),
\item for each $v\in V$, the disks $D_e^-$ for $e\in E_v^-$
and $D_e^+$ for $e\in E_v^+$ are disjoint in $S_v$,
\item there exists positive constants $\cc$ and $\ccc$ such that
$$\forall e\in E,\quad\forall z\in D_e^{\pm},
\quad \cc|z_e^{\pm}(z)|\leq |z-p_e^{\pm}|\leq \ccc|z_e^{\pm}(z)|.$$
\end{enumerate}
\end{definition}
\begin{remark}
\label{remark-admissible}
If the coordinates are admissibles, then for each edge $e\in E$,
the round disk $D(p_e^{\pm},\rho\cc)$ is included in the topological
disk $D_e^{\pm}$.
In particular, for each vertex $v\in V$, the disks $D(p_e^+,\rho\cc)$
for $e\in E_v^+$ and $D(p_e^-,\rho\cc)$ for $e\in E_v^-$ are disjoint.
\end{remark}
\subsection{The $\ell^{\infty}$ case}
Next, as $\Sigma_{\tt}$ is non compact, the space of regular differentials
is infinite
dimensional so we need to define a norm on this space.
Fix some real number $0<\epsilon\leq \rho\cc$.
For each vertex $v\in V$, let $\Omega_{v,\epsilon}$ be the Riemann sphere $S_v$ minus the
disks $D(p_e^-,\epsilon)$ for $e\in E_v^-$ and $D(p_e^+,\epsilon)$ for
$e\in E_v^+$.
Let $\boH^1_{\infty}(\Sigma_{\tt})$
be the space of regular differentials
on $\Sigma_{\tt}$ such that the norm
$$||\omega||_{\infty}
=\sup_{v\in V}\sup_{\Omega_{v,\epsilon}}\left|\frac{\omega}{dz}\right|$$
is finite.
Let $\ell^{\infty}(E)$ be the space of bounded sequences of complex numbers
$(\alpha_e)_{e\in E}$ with the $\sup$ norm.
Our first result, which generalises theorem \ref{theorem-fay} to the infinite case, is
\begin{theorem}
\label{theorem1}
Assume that the degree of the vertices of $\Gamma$ is bounded and
that the coordinates are admissible.
Then for $\tt$ small enough (in $\ell^{\infty}$ norm), 
the operator
$\omega\to (\int_{\gamma_e}\omega)_{e\in E}$
is an isomorphism of Banach spaces
from $\boH^1_\infty(\Sigma_{\tt})$ to the subspace of sequences in
$\ell^{\infty}(E)$ which satisfy \eqref{eq-compatible}.
\end{theorem}
\begin{remark}
\label{remark1}
If follows from the theorem that the space
$\boH^1_{\infty}(\Sigma_{\tt})$
does not depend on the choice of $\epsilon$ (taking another $\epsilon$ will define
an equivalent norm).
\end{remark}
\subsection{Examples of admissible coordinates}
\begin{example}
\label{example1}
\em Assume that for each $v\in V$,
the points $p_e^+$ for $e\in E_v^+$ and $p_e^-$ for $e\in E_v^-$
are at distance greater than $2\rho$ from each other.
We can take the coordinates
$z_e^{\pm}(z)=z-p_e^{\pm}$ on $D_e^{\pm}=D(p_e^{\pm},\rho)$.
These coordinates are admissible with $\cc=\ccc=1$.
\end{example}

This is probably the most natural choice of coordinates. The proof
of theorem \ref{theorem2} would be significantly simplified if we restricted
ourselves to these coordinates, so let me give other examples
to motivate the use of more general coordinates.

\begin{example}
\label{example2}
\em Assume that for each $v\in V$, the points $p_e^+$ for $e\in E_v^+$
and $p_e^-$ for $e\in E_v^-$
are on the unit circle and are at distance greater than $4 \rho$
from each other.
We would like the inversion $\sigma(z)=1/\overline{z}$ to be
well defined on $\Sigma_{\tt}$.
In this case we take the coordinates
$$z_e^{\pm}(z)=i\;\frac{z-p_e^{\pm}}{z+p_e^{\pm}}.$$
Then $z_e^{\pm}\circ\sigma=\overline{z_e^{\pm}}$, so provided all $t_e$
are real numbers, $\sigma$ is well defined on $\Sigma_{\tt}$.
Let $D_e^{\pm}$ be the inverse image of $D(0,\rho)$ by
$z_e^{\pm}$.
Computations shows that these disks
are disjoint, and the coordinates are admissible with
$\cc=\frac{2}{1+\rho}$ and
$\ccc=\frac{2}{1-\rho}$.
\end{example}
\begin{example}
\label{example3}
\em
Assume that the graph $\Gamma$ is bipartite, meaning that the set of
vertices has a partition $V=V^+\cup V^-$, so that all edges go from one
vertex in $V^-$ to a vertex in $V^+$.
Consider on each sphere $S_v$ a meromorphic function $f_v$ with simple
zeros at $p_e^-$ for $e\in E_v^-$ and $p_e^+$ for $e\in E_v^+$.
For each $e\in E$ from 
$v$ to $v'$, we take $z_e^-=f_v$ and $z_e^+=f_{v'}$.
Some hypothesis must be made on the functions $(f_v)_{v\in V}$
so that these coordinates are admissibles.
Consider a small 
complex number $t$ and take $t_e=t^2$ for all $e\in E$. Then we
can define a meromorphic function $f$ on $\Sigma_{\tt}$ by
$f(z)=\disp\frac{f_v(z)}{t} \mbox{ if } v\in V^+$
and $f(z)=\disp \frac{t}{f_v(z)} \mbox{ if } v\in V^-$.
Indeed, for each edge $e\in E$ from $v$ to $v'$, if $z\in D_e^-$ is
identified with $z'\in D_e^+$, we have $f_v(z) f_{v'}(z') = t^2$ so
$f(z)=f(z')$ and $f$ is well defined.

\medskip

The reason this example is interesting is that together with a Riemann surface, we construct a meromorphic function.
It is in general not easy to define a meromorphic function
on a given Riemann surface (think of Abel's theorem).
This is particularly fruitful in the case of minimal surfaces, where the
function $f$ will be the Gauss map. See \cite{traizet} for an illustration of
this idea.
\end{example}
\subsection{Admissible norms}
\label{section-admissiblenorms}
The $\ell^{\infty}$ norm is maybe the most natural one, but for certain applications
it is desirable to allow other norms. For example, in section \ref{section-decay}
we use weighted $\ell^{\infty}$ norms to study the decay of regular differentials.
In the following definition, we try to consider norms as general as possible.

Let $\C^E$ and $\C^V$ denote the space of sequences of complex numbers
$(\alpha_e)_{e\in E}$ and $(u_v)_{v\in V}$ indexed by edges and vertices,
respectively.
We assume that we are given two norms $||\cdot||_E$ and $||\cdot||_V$ defined
on a subspace $\boB_E$ of $\C^E$ and $\boB_V$ of $\C^V$.
\begin{definition}
\label{def-admissiblenorms}
We say that the norms $||\cdot||_E$ and $||\cdot||_V$ are
admissible if the following
conditions hold :
\begin{enumerate}
\item The spaces $\boB_E$ and $\boB_V$ are Banach spaces.
\item The norms $||\cdot||_E$ and $||\cdot||_V$ are monotonic, in the
sense that if $|\alpha_e|\leq |\alpha'_e|$ for all $e\in E$ then
$||\alpha||_E\leq ||\alpha'||_E$, and a similar definition for $||\cdot||_V$.
\item The operator
$\disp (\alpha_e)_{e\in E}\mapsto \left(\sum_{e\in E_v} |\alpha_e|\right)_{v\in V}$
is bounded from $\boB_E$ to $\boB_V$.
\item The operator
$\disp (u_v)_{v\in V}\mapsto \left(|u_v|+|u_{v'}|
\right)_{e\in E\atop e=vv'}$
is bounded from $\boB_V$ to $\boB_E$.
\item The space of sequences $(\lambda_{e,n})_{e\in E,n\in\N}$ for
which the norm $\disp ||(\sup_{n\in \N} |\lambda_{e,n}|)_{e\in E}||_{E}$
is defined is a Banach space.
\end{enumerate}
\end{definition}
We define $\boH^1(\Sigma_{\tt})$ as the space of regular differentials on
$\Sigma_{\tt}$ such that the norm
$$||\omega||_{\Omega}=
\left\|\left(\sup_{\Omega_{v,\epsilon}}\left|\frac{\omega}{dz}\right|\right)_{v\in V}\right\|_{V}
$$ is defined.
The following theorem, which is the main result of this paper,
generalises theorem \ref{theorem1} to this more general setup.
\begin{theorem}
\label{theorem2}
Assume that the coordinates and the norms are admissible.
Then for 
$\tt$ small enough (in $\ell^{\infty}$ norm),
the operator $\omega\mapsto (\int_{\gamma_e}\omega)_{e\in E}$ is an isomorphism
from $\boH^1(\Sigma_{\tt})$ to the subspace of
sequences in $\boB_E$ which satisfy \eqref{eq-compatible}
\end{theorem}
\subsection{Examples of admissible norms}
We define weighted $\ell^p$ spaces for $1\leq p\leq \infty$ as follows.
Consider a function $\sigma:V\to (0,\infty)$.
Given a sequence $u=(u_v)_{v\in V}$, let
$$||u||_{p,\sigma}=\left(\sum_{v\in V}\sigma(v)^p |u_v|^p\right)^{1/p}\qquad
\mbox{ if } p<\infty,$$
$$||u||_{\infty,\sigma}=\sup_{v\in V}\sigma(v) |u_v|.$$
We call $\ell^{p,\sigma}(V)$ the space of sequences $u\in \C^V$ for which the above norm
is finite.
When $\sigma\equiv 1$, these are the usual $\ell^p$ spaces.

We define a weight function on edges by
$\sigma(e)=\frac{1}{2}(\sigma(v)+\sigma(v'))$ if $e$ is an edge from $v$ to $v'$.
The space $\ell^{p,\sigma}(E)$ of sequences $\alpha\in\C^E$ which have finite
$\ell^{p,\sigma}$ norm is defined in the obvious way.

\begin{proposition}
\label{proposition-lpnorms}
Assume that there exists $k$ such that the degree of the vertices of $\Gamma$ is bounded by $k$,
and that for any adjacent vertices $v$ and $v'$, $\sigma(v)\leq k \sigma(v')$ (and
$\sigma(v')\leq k\sigma(v)$). Then the norms $||\cdot||_E=||\cdot||_{p,\sigma}$ and
$||\cdot||_V=||\cdot||_{p,\sigma}$ are admissible.
\end{proposition}
Proof.
The first point is standard, the second is clear. Regarding the third point, in the case
$p<\infty$ we write
\begin{eqnarray*}
\sum_{v\in V}\sigma(v)^p\left(\sum_{e\in E_v} |\alpha_e|\right)^p
&\leq&\sum_{v\in V} \sigma(v)^p \deg(v)^{p-1}\sum_{e\in E_v} |\alpha_e|^p
\quad \mbox{ by Jensen inequality}\\
&\leq&  k^{p-1} \sum_{v\in V} \sum_{e\in E_v} \sigma(v)^p |\alpha_e|^p\\
&=&k^{p-1}\sum_{e\in E} (\sigma(v)^p+\sigma(v')^p) |\alpha_e|^p
\\
&\leq& 2^p k^{p-1}\sum_{e\in E} \sigma(e)^p |\alpha_e|^p
\end{eqnarray*}
The case $p=\infty$ is straightforward.

Regarding the fourth point, we have for any edge $e\in E$ from $v$ to $v'$,
$\sigma(e)\leq \frac{k+1}{2}\sigma(v)$ and
$\sigma(e)\leq \frac{k+1}{2}\sigma(v')$, so
\begin{eqnarray*}
\sum_{e\in E\atop e=vv'} \sigma(e)^p (|u_v|+|u_{v'}|)^p
&\leq& \sum_{e\in E} \sigma(e)^p 2^{p-1} (|u_v|^p + |u_{v'}|^p)\\
&\leq& \frac{(1+k)^p}{2}\sum_{e\in E} \sigma(v)^p |u_v|^p+\sigma(v')^p |u_{v'}|^p\\
&=& \frac{(1+k)^p}{2}\sum_{v\in V} \deg(v) \sigma(v)^p |u_v|^p\\
&\leq & \frac{k(1+k)^p}{2} \sum_{v\in V} \sigma(v)^p |u_v|^p
\end{eqnarray*}

Finally, point 5 is true by the standard fact that $\ell^{\infty}(\N)$ is a Banach space and
an infinite product of Banach spaces is a Banach
space for the $\ell^p$ norm of the norms
(cf exercice T page 243 in \cite{brown-pearcy}).
\cqfd
\begin{remark}
Theorem \ref{theorem1} is a corollary of theorem \ref{theorem2}
and proposition \ref{proposition-lpnorms}.
\end{remark}
\subsection{Decay at infinity of normalised differentials}
\label{section-decay}
To illustrate the use of weights, let us consider the following example.
\begin{example}
\label{example4}
\em
Let $\Gamma$ be an infinite graph which contains at least one cycle
$\gamma$.
Give the cycle $\gamma$ an orientation and orient the rest of $\Gamma$ arbitrarily.
Define the sequence $(\alpha_e)_{e\in E}$ as $\alpha_e=1$ if $e$ belongs
to $\gamma$ and $\alpha_e=0$ otherwise.
Then $(\alpha_e)_{e\in E}$ satisfies \eqref{eq-compatible}, so defines
a regular differential $\omega\in\boH^1_{\infty}(\Sigma_{\tt})$. We call
$\omega$ the normalised differential associated to the cycle $\gamma$.  
\end{example}
When $\tt=0$, $\omega$ is supported on the spheres $S_v$ such that
the vertex $v$ belongs to the cycle $C$, so has compact support. It is natural to ask 
what is the decay of $\omega$ at infinity when $\tt\neq 0$.
For this purpose, we pick an arbitrary vertex $v_0\in V$ and consider the weight
$\sigma(v)=r^{d(v,v_0)}$, where $r>1$ is a fixed real number, and $d(v,v_0)$
denotes the graphical distance from $v_0$ to $v$ on $\Gamma$.
Then $\alpha\in \ell^{\infty,\sigma}(E)$, so by the theorem,
$\omega\in \boH^1_{\infty,\sigma}(\Sigma_{\tt})$ for $\tt$ small enough,
so $\omega$ has exponential decay.
More precisely $\sup_{\Omega_{v,\epsilon}}\left|\frac{\omega}{dz}\right|$ is bounded by some constant times $r^{-d(v,v_0)}$.
The rate of decay $r$ is arbitrary, but of course the larger $r$, the smaller
$\tt$ must be.
\section{Proof of theorem 3}
\label{section-proof}
\begin{convention}
By a uniform constant, we mean a number which only depends on the numbers
$\rho$,
$\cc$ and $\ccc$ in definition \ref{def-admissiblecoordinates} and
the norm of the operators in points 3 and 4 of definition
\ref{def-admissiblenorms}.
We use the letter $C$ to denote uniform constants.
The same letter $C$ can be used to denote various uniform constants.
\end{convention}

\medskip

Let $L$ be the operator $L(\omega)=(\int_{\gamma_e}\omega)_{e\in E}$.
We prove that
\begin{enumerate}
\item $L$ is bounded,
\item $L$ is onto,
\item $L$ is injective.
\end{enumerate}

1) Let $e$ be an edge from $v$ to $v'$.
Let $r=\rho\cc\geq\epsilon$.
By remark \ref{remark-admissible},
we can replace
$\gamma_e$ by the circle $C(p_e^+,r)$, so
$$\left|\int_{\gamma_e}\omega\right|=\left|\int_{C(p_e^+,r)}
\omega\right|
\leq 2\pi r\sup_{\Omega_{v',\epsilon}}|\frac{\omega}{dz}|
\leq 2\pi r\left(
\sup_{\Omega_{v,\epsilon}}|\frac{\omega}{dz}|
+
\sup_{\Omega_{v',\epsilon}}|\frac{\omega}{dz}|
\right).$$
By point 4 of the definition of admissible norms, we get
$||L(\omega)||_E\leq  C||\omega||_{\Omega}$ for some uniform
constant $C$.
\medskip

2) Let us prove next that $L$ is onto.
Given $(\alpha_e)_{e\in E}\in \boB_E$, satisfying \eqref{eq-compatible},
we are asked to construct a regular differential $\omega$ on $\Sigma_{\tt}$
such that $L(\omega)=\alpha$.
In each sphere $S_v$, we define $\omega=\omega_1+\omega_2$ with
\begin{equation}
\label{eq-omega1}
\omega_1=
\omega_1^+ + \omega_1^-
=
\frac{1}{2\pi i}\sum_{e\in E_v^+} \alpha_e\frac{dz}{z-p_e^+} 
-\frac{1}{2\pi i}\sum_{e\in E_v^-} \alpha_e\frac{dz}{z-p_e^-} 
\end{equation}
\begin{equation}
\label{eq-omega2}
\omega_2=
\omega_2^+ + \omega_2^-
=
\sum_{e\in E_v^+} \sum_{n=2}^{\infty} \lambda_{e,n}^+ 
\frac{
\left(\frac{\epsilon}{4}\right)^n
dz}{(z-p_e^+)^n}
+\sum_{e\in E_v^-} \sum_{n=2}^{\infty} \lambda_{e,n}^-
\frac{
\left(\frac{\epsilon}{4}\right)^n
dz}{(z-p_e^-)^n}.
\end{equation}
($\omega_1^+$ and $\omega_1^-$ denote the first and second sum, respectively.)
The meromorphic differential $\omega_1$ has the required periods and is
entirely determined by the data we are given.
Condition \eqref{eq-compatible} ensures that the residue of $\omega_1$
at infinity vanishes, so $\omega_1$ is holomorphic at infinity.
The meromorphic differential $\omega_2$ is a corrective term with no periods.
The coefficients $\lambda_{e,n}^{\pm}$ are to be determined.

\medskip
Let us first assume that the series defining $\omega_2$ converge
and see what is the condition so that $\omega$ is a well defined
differential on $\Sigma_{\tt}$.
Then we solve theses equations, and finally prove that the series
converge.

\medskip

For any edge $e\in E$, define
$$\varphi_e^+=(z_e^+)^{-1}\circ\left(\frac{t_e}{z_e^-}\right)
,\qquad
\varphi_e^-=(z_e^-)^{-1}\circ\left(\frac{t_e}{z_e^+}\right)
=(\varphi_e^+)^{-1}$$
If $t_e\neq 0$,
$\varphi_e^+$ is a biholomorphism from the annulus
$\frac{|t_e|}{\rho}\leq |z_e^-|\leq \rho$ to the annulus $\frac{|t_e|}{\rho}\leq |z_e^+|\leq \rho$,
and the point $z$ is identified with the point
$\varphi_e^+(z)$ when defining $\Sigma_{\tt}$.
In other words, these are the change of charts for the coordinates
$z_e^+$ and $z_e^-$ on $\Sigma_{\tt}$,
so $\omega$ is well defined on $\Sigma_{\tt}$ if
$(\varphi_e^+)^*\omega=\omega$ for all edges $e\in E$.
\begin{claim}
\label{claim1}
$\omega$ is well defined on $\Sigma_{\tt}$ if and only if
for all edges $e\in E$ such that $t_e\neq 0$
$$\forall n\geq 0,\quad
\int_{\partial D_e^-} (z-p_e^-)^n( (\varphi_e^+)^*\omega-\omega)=0,$$
$$\forall n\geq 0,\quad
\int_{\partial D_e^+} (z-p_e^+)^n( (\varphi_e^-)^*\omega-\omega)=0.$$
\end{claim}
Proof. By the theorem on Laurent series,
the holomorphic differential $(\varphi_e^+)^*\omega-\omega$ is zero
on the annulus $\frac{|t_e|}{\rho}\leq |z_e^-|\leq \rho$ if and only if for
all $n\in\Z$,
\begin{equation}
\label{eq-welldefined}
\int_{\partial D_e^-}(z_e^-)^n ((\varphi_e^+)^*\omega-\omega)=0.
\end{equation}
Equation \eqref{eq-welldefined} for all $n\geq 0$ means that
$(\varphi_e^+)^*\omega-\omega$ extends holomorphically to the
disk $|z_e^-|\leq \rho$. Using $z-p_e^-$ as a coordinate, this is equivalent
to the first condition of the claim.
If $n\leq 0$, then by a change of variable
\begin{eqnarray*}
\int_{\partial D_e^-}(z_e^-)^n ((\varphi_e^+)^*\omega-\omega)
&=& 
-\int_{\partial D_e^+}(\varphi_e^-)^*
\left[(z_e^-)^n ((\varphi_e^+)^*\omega-\omega)\right]
\\
&=&
-\int_{\partial D_e^+}
\left(\frac{t_e}{z_e^+}\right)^n(\omega-(\varphi_e^-)^*\omega)
\end{eqnarray*}
So equation \eqref{eq-welldefined} for all $n\leq 0$ means that
$\omega-(\varphi_e^-)^*\omega$ extends holomorphically to the disk
$|z_e^+|\leq \rho$. Using $z-p_e^+$ as a coordinate, this is equivalent to
the second condition of the claim.
\cqfd
\begin{claim}
\label{claim2}
$\omega$ is a well defined regular differential 
on $\Sigma_{\tt}$ if and only
if for all $e\in E$ and all $n\geq 2$,
\begin{equation}
\label{eq-lambda-}
\lambda_{e,n}^-=\frac{-1} {2\pi i}
\left(\frac{4}{\epsilon}\right)^n
\int_{\partial D_e^+}(\varphi_e^- -p_e^-)^{n-1}\omega,
\end{equation}
\begin{equation}
\label{eq-lambda+}
\lambda_{e,n}^+=\frac{-1}
{2\pi i}
\left(\frac{4}{\epsilon}\right)^n
\int_{\partial D_e^-}(\varphi_e^+-p_e^+)^{n-1}\omega.
\end{equation}
\end{claim}
Proof.
If $t_e=0$, then $\varphi_e^+-p_e^+$ and $\varphi_e^- -p_e^-$ are 
identically zero.
Hence equations \eqref{eq-lambda-} and \eqref{eq-lambda+} for all
$n\geq 2$ mean that $\omega_2$ is holomorphic at $p_e^-$ and 
$p_e^+$, so $\omega$ has at most simple poles as required for a
regular differential.

If $t_e\neq 0$, we use the previous claim.
The definition of $\omega=\omega_1+\omega_2$ gives
$$\int_{\partial D_e^-}(z-p_e^-)^n\omega
=\left\{\begin{array}{l}
-\alpha_e \mbox{ if } n=0\\
2\pi i \lambda_{e,n+1}^- \left(\frac{\epsilon}{4}\right)^{n+1}
\mbox{ if } n\geq 1
\end{array}\right.$$
By a change of variable,
\begin{eqnarray*}
\int_{\partial D_e^-} (z-p_e^-)^n (\varphi_e^+)^*\omega
&=&
-\int_{\partial D_e^+} (\varphi_e^-)^* \left[
(z-p_e^-)^n (\varphi_e^+)^*\omega\right]\\
&=&
-\int_{\partial D_e^+} (\varphi_e^- -p_e^-)^n\omega\\
&=& -\alpha_e \mbox{ if } n=0
\end{eqnarray*}
We have similar statements exchanging the roles of the $+$ and $-$ signs.
The claim follows.
\cqfd

Given a holomorphic differential $w$ on $\Omega_{\epsilon}=\bigcup_{v\in V} \Omega_{v,\epsilon}$,
define for $e\in E$ and $n\geq 2$
\begin{equation}
\label{eq-F}
F_{e,n}^{\pm}(w)=\frac{-1}{2\pi i}\left(\frac{4}{\epsilon}\right)^n\int_{\partial D_e^{\mp}}
(\varphi_e^{\pm}-p_e^{\pm})^{n-1}w.
\end{equation}
Let $F(w)=(F_{e,n}^{\pm}(w))_{e\in E,n\geq 2}$.
By claim \ref{claim2}, $\omega$ is well defined on $\Sigma_{\tt}$ if
$\lambda=F(\omega)$.
We want to find $\lambda$ as a fixed point of the map
$\lambda\mapsto F(\omega_1+\omega_2(\lambda))$.
To this effect, let $\boB_L$
be the space of sequences $\lambda=(\lambda_{e,n}^{\pm})_{e\in E,n\geq 2}$
such that the following norm is defined :
$$||\lambda||_L=\left\|\left(\sup_{n\geq 2} \max\{
|\lambda_{e,n}^+|,|\lambda_{e,n}^-|\}\right)_{e\in E}\right\|_E.$$
This is a Banach space by point 5 of definition \ref{def-admissiblenorms}.
\begin{lemma}
\label{lemma1}
There exists uniform constants such that if
$||\tt||_{\infty}\leq \frac{\cc\rho\epsilon}{4\ccc^2}$,
\begin{enumerate}
\item $||\omega_1||_{\Omega}\leq \frac{C}{\epsilon} ||\alpha||_E$,
\item $||\omega_2||_{\Omega}\leq C ||\lambda||_L$,
\item $||F(w)||_L\leq \frac{C}{\epsilon^2} ||\tt||_{\infty}||w||_{\Omega}$.
\end{enumerate}
\end{lemma}
Proof of point 1 :
we have the straightforward estimate
$$\sup_{\Omega_{v,\epsilon}}\left|\frac{\omega_1}{dz}\right|
\leq \frac{1}{\epsilon}\sum_{e\in E_v}|\alpha_e|.$$
The conclusion follows by point 3 of definition \ref{def-admissiblenorms}.

\medskip

\noindent Proof of point 2 :
$$\sup_{\Omega_{v,\epsilon}}\left|\frac{\omega_2^+}{dz}\right|
\leq
\sum_{e\in E_v^+}\sum_{n=2}^{\infty}
\frac{|\lambda_{e,n}^+|}{\epsilon^n}
\left(\frac{\epsilon}{4}\right)^n
\leq
\left(\sum_{n=2}^{\infty}\left(\frac{1}{4}\right)^n\right)
\left(\sum_{e\in E_v^+}
\sup_{n\geq 2}|\lambda_{e,n}^+|\right). $$
We estimate $\omega_2^-$ in the same way.
The conclusion follows by point 3 of definition \ref{def-admissiblenorms}.

\medskip

\noindent Proof of point 3 :
Let $e$ be an edge from $v$ to $v'$.
By definition of admissible coordinates, we have
$$|\varphi_e^+(z)-p_e^+|
\leq \ccc |z_e^+(\varphi_e^+(z))|
=\frac{\ccc|t_e|}{|z_e^-(z)|}
\leq \frac{\ccc^2|t_e|}
{|z-p_e^-|}.$$
Arguing in the same way on the other side we obtain the following
useful estimate
\begin{equation}
\label{eq-varphi}
\frac{\cc^2 |t_e|}{|z-p_e^-|}
\leq |\varphi_e^+(z)-p_e^+|\leq
\frac{\ccc^2 |t_e|}{|z-p_e^-|}
\end{equation}
We replace the circle
$\partial D_e^-$ in the definition of $F_{e,n}^+$ by
the circle $C(p_e^-,r)$ with $r=\cc\rho$.
Then we have for $e\in E$ and $n\geq 2$
\begin{eqnarray}
\label{eq-estimee-F}
|F_{e,n}^+(w)|
&\leq & \frac{1}{2\pi}2\pi r \frac{4}{\epsilon}\left(\frac{4}
{\epsilon}\frac{\ccc^2 |t_e|}{r}
\right)^{n-1}
\sup_{\Omega_{v,\epsilon}}\left|\frac{w}{dz}\right|
\\
\nonumber
&\leq &  \frac{4r}{\epsilon}
\left(
\frac{4\ccc^2 |t_e|}{\epsilon r}
\right)
\sup_{\Omega_{v,\epsilon}}\left|\frac{w}{dz}\right|
\\
\nonumber
&\leq & \frac{16\ccc^2}{\epsilon^2}
||\tt||_{\infty} \sup_{\Omega_{v,\epsilon}}\left|\frac{w}{dz}\right|.
\end{eqnarray}
(The term in the middle parenthesis of the first line
is less than one by our hypothesis
on $\tt$).
We estimate $F_{e,n}^-(w)$ in the same way.
The conclusion follows by point 4 of definition \ref{def-admissiblenorms}.
\cqfd

It follows from the lemma that if $||\tt||_{\infty}$ is small enough
(depending on $\epsilon$),
the operator $\lambda\mapsto F(\omega_1+\omega_2(\lambda))$
is contracting from $\boB_L$ to itself. By the contraction
mapping theorem, there exists a unique $\lambda\in \boB_L$
such that $F(\lambda)=\lambda$.

\medskip

It remains to consider the convergence of the series defining
$\omega_2$.
First of all, we have actually seen in the proof of point 2 of the lemma that
this series converges on $\Omega_{v,\epsilon}$, so $\omega$ is already defined on each $\Omega_{v,\epsilon}$.
This is not enough because these domains do not cover all of $\Sigma_{\tt}$.
Since $||\omega||_{\Omega}$ is finite, $\sup_{\Omega_{v,\epsilon}}|\omega|$ is finite for
all $v\in V$ (altough maybe not uniformly bounded).
Since $\lambda=F(\omega)$, we have by equation \ref{eq-estimee-F}
$$
\left(\frac{\epsilon}{4}\right)^n|\lambda_{e,n}^{\pm}|\leq r\left(\frac{\ccc^2 |t_e|}
r\right)^{n-1}\sup_{\Omega_{v,\epsilon}}\left|\frac{\omega}{dz}\right|.$$
It follows that the series defining $\omega_2$ in $S_v$ converges provided
$|z-p_e^{\pm}|>\frac{\ccc^2 |t_e|}{r}$. Since these domains
cover all of $\Sigma_{\tt}$ for $\tt$ small enough,
$\omega$ is a well defined regular differential on $\Sigma_{\tt}$.
So we have proven that $L$ is onto.
\begin{remark}
\label{remark2}
It follows from the lemma that for $||\tt||_{\infty}$
small enough (depending on $\epsilon$),
$||\omega||_{\Omega}\leq \frac{C}{\epsilon} ||\alpha||_{E}$.
\end{remark}

\medskip

\noindent 3) Finally we prove that $L$ is injective for $||\tt||_{\infty}$
small enough.
Let $\omega$ be in the kernel of $L$.
For any $e\in E$, $\omega$ is holomorphic in the annulus
$\ccc\frac{|t_e|}{\rho}\leq |z-p_e^{\pm}|\leq \cc\rho$
(which is included in the annulus
$\frac{|t_e|}{\rho}\leq |z_e^{\pm}|\leq \rho$),
so we can write its Laurent series in
this annulus as
$$\omega=\sum_{n=-\infty}^{\infty} c_{e,n}^{\pm} \frac{dz}{(z-p_e^{\pm})^n}.$$
The coefficients $c_{e,1}^{\pm}$ are all zero because $L(\omega)=0$.
The sum for $n\leq 0$ extends holomorphically to the disk
$D(p_e^{\pm},r)$.
The sum for $n\geq 2$ extends holomorphically to the outside of this disk.
Therefore, the difference
$$\omega-
\sum_{e\in E_v^+}
\sum_{n=2}^{\infty} c_{e,n}^+ \frac{dz}{(z-p_e^+)^n}
-
\sum_{e\in E_v^-}
\sum_{n=2}^{\infty} c_{e,n}^- \frac{dz}{(z-p_e^-)^n}$$
extends holomorphically to all of $S_v$. Since there are no holomorphic
1-form on the sphere, it is zero.
Let $\lambda_{e,n}^{\pm}=(\frac{4}{\epsilon})^n c_{e,n}^{\pm}$, then $\omega=\omega_2(\lambda)$.
Since $\omega$ is well defined on $\Sigma_{\tt}$, we
have $\lambda=F(\omega)$.
By lemma \ref{lemma1},
$||\lambda||_L\leq \frac{C}{\epsilon^2} ||\tt||_{\infty} ||\lambda||_L$,
so $\lambda=0$
if $\tt$ is small enough. Hence $\omega=0$.
This concludes the proof of theorem \ref{theorem2}.
\cqfd
\begin{remark}
\label{remark3}
We have proven in point 2 that any regular differential
$\omega\in\boH^1(\Sigma_{\tt})$ can be written $\omega=
\omega_1(\alpha)+\omega_2(\lambda)$,
where $\omega_1$ and $\omega_2$ are defined by \eqref{eq-omega1} and
\eqref{eq-omega2}. Moreover, $\lambda=F(\omega)$, where $F(\omega)$
is defined by \eqref{eq-F}.
\end{remark}
\section{Smooth dependance on parameters}
\label{section-smooth}
By theorem \ref{theorem2}, for any $\alpha\in\boB_E$ satisfying \eqref{eq-compatible},
there exists a unique regular differential $\omega$ on $\Sigma_{\tt}$ whose periods are prescribed
by $\alpha$. In this section, we prove that $\omega$ depends smoothly on
the parameters in the constuction of $\Sigma_{\tt}$, namely
$\pp=(p_e^{\pm})_{e\in E}$ and $\tt=(t_e)_{e\in E}$.
Before we can formulate this result, several points need to be adressed.
\medskip

We assume that the parameter $\pp=(p_e^{\pm})_{e\in E}$ is
in a neighborhood of some central value
$\underline{\pp}=(\underline{p}_e^{\pm})_{e\in E}$
in $\ell^{\infty}$ norm.
In this context, we define $\Omega_{v,\epsilon}$ to be the sphere $S_v$
minus the disks $D(\underline{p}_e^-,\epsilon)$ for $e\in E_v^-$ and
$D(\underline{p}_e^+,\epsilon)$ for $e\in E_v^+$.
We define $\Omega_{\epsilon}$ as the disjoint union of all $\Omega_{v,\epsilon}$ for $v\in V$.
The norm $||\cdot||_{\Omega}$ is defined as in section \ref{section-admissiblenorms}.
The point here is that the domain $\Omega_{\epsilon}$ is fixed,
whereas our former domain $\Omega_{\epsilon}$ depends on the
parameter $\pp$.
Thanks to remark \ref{remark1}, the new norm $||\cdot||_{\Omega}$
on $\boH_1(\Sigma_{\tt})$
is equivalent to the former one.

\medskip

The coordinates $z_e^{\pm}$ must depend in some ways on the parameters, else
the requirement $z_e^{\pm}(p_e^{\pm})=0$ fixes the point $p_e^{\pm}$.
We assume that for each edges $e\in E$, the coordinates $z_e^{\pm}$ have the form
$z_e^{\pm}(z)=\zeta_e^{\pm}(z-p_e^{\pm},\xi_e^{\pm})$, where $\xi_e^{\pm}$ is a vector of
complex parameters and $\zeta_e^{\pm}$ is a holomorphic function in all its variables.
For example, in example \ref{example2}, we would take
$\xi_e^{\pm}=p_e^{\pm}$ and
$\zeta_e^{\pm}(z,\xi)=i\frac{z}{z+2\xi}$.

\medskip

We assume that each parameter $\xi_e^{\pm}$ is in a finite
dimensional complex vector space, possibly depending on the edge,
but of uniformly bounded dimension.
We write $\xi=(\xi_e^{\pm})_{e\in E}$ and we assume that
the coordinates $z_e^{\pm}$ are defined and admissible
for all values of $\xi$ in a neighborhood of some
central value $\underline{\xi}=(\underline{\xi}_e^{\pm})_{e\in E}$
in $\ell^{\infty}$ norm.

\begin{theorem}
\label{theorem-smooth}
The regular differential $\omega$, restricted to $\Omega_{\epsilon}$,
depends smoothly on the parameters
$(\pp,\tt,\xi)$ in a neighborhood of $(\underline{\pp},0,\underline{\xi})$ in $\ell^{\infty}$
norm.
\end{theorem}
The norm on the domain space is the $\ell^{\infty}$ norm. The norm on the target
space is the norm $||\cdot||_{\Omega}$.
\begin{remark}
\label{remark4}
By a theorem of Graves-Taylor-Hille-Zorn, a map
$f$ from an open set of a {\em complex} Banach space $E$ to a complex
Banach space $F$, which is differentiable in the usual (Frechet) sense,
is holomorphic, hence smooth.
So the word ``smoothly'' in the theorem can be replaced by the
word ``analytically''.
See the book \cite{soo} for the definition of a holomorphic map
between complex Banach spaces and theorem 14.3 for the statement and
the proof.
Also, Hartog's theorem holds in the complex Banach space setup :
if $E_1$ and $E_2$ are Banach spaces and $f:E_1\times E_2\to F$ is
separately holomorphic with respect to each of its two variables,
then it is holomorphic (theorem 14.27).
\end{remark}

\medskip

\noindent
Proof of theorem \ref{theorem-smooth}
(continued from the proof of theorem \ref{theorem2}) :
recall that we found $\lambda$ as a fixed point of $\lambda\mapsto F(\omega_1)
+F(\omega_2(\lambda))$. By the following lemma and the fixed point theorem with parameters, $\lambda$
depends smoothly on all parameters,
so $\omega=\omega_1+\omega_2(\lambda)$
also depends smoothly on all parameters.
\cqfd

\medskip

Let $\boH_1(\Omega_{\epsilon})$ be the space of holomorphic differentials $w$
on $\Omega_{\epsilon}$ such that $||w||_{\Omega}$ is defined.
\begin{lemma}
\label{lemma-smooth}
For $(\pp,\tt,\xi)$ in a neighborhood of $(\underline{\pp},0,\underline{\xi})$
in $\ell^{\infty}$ norm, $\alpha\in\boB_E$, $\lambda\in\boB_L$
and $w\in\boH_1(\Omega_{\epsilon})$,
\begin{enumerate}
\item the map $(\pp,\alpha)\mapsto \omega_1(\pp,\alpha)$ is smooth,
\item the map $(\pp,\lambda)\mapsto \omega_2(\pp,\lambda)$ is smooth,
\item the map $(\pp,\tt,\xi,w)\mapsto F(\pp,\tt,\xi,w)$ is smooth.
\end{enumerate}
\end{lemma}
To prove this lemma, we first consider an abstract result.
The setup is the following : we have an infinite family
of holomorphic functions $f_n(x_n)$ for $n\in\N$.
The variable $x_n$ is in the polydisk $B(0,R)$ in $\C^{d_n}$.
We assume that the dimension $d_n$ is uniformly bounded.
The function $f_n$ takes value in a Banach space $F_n$.
Let $E=\prod_{n\in\N} \C^{d_n}$ and $F=\prod_{n\in\N} F_n$, both with the sup norm.
Let $x=(x_n)_{n\in\N}$ and $f(x)=(f_n(x_n))_{n\in\N}$
\begin{lemma}
\label{lemma3}
If each $||f_n||$ is bounded on $B(0,R)$
by a constant $C$ independant of $n$, then
$f:B(0,R)\subset E\to F$ is smooth.
\end{lemma}
Proof. Consider some $r<R$. Since the function $f_n$ is holomorphic
in the polydisk $B(0,R)$ and bounded by $C$,
all its partial derivatives of order $k$ in the polydisk $B(0,r)$
have norm bounded by $C k! \,(R-r)^{-k}$ by Cauchy's estimates, which holds
true for Banach valued holomorphic functions.
The point is that this bound does not depend on $n$.
It is then straightforward to check that $f$ is differentiable
in the ball $B(0,r)$, with differential $df(x)(h)=(df_n(x_n)(h_n))_{n\in\N}$.
Since $r$ is arbitrary, $f$ is differentiable in the open ball $B(0,R)$.
The fact that $f$ is smooth follows by induction (or using
the theorem of Graves-Taylor-Hille-Zorn.)
\cqfd

\medskip

\noindent
Proof of point 1 of lemma \ref{lemma-smooth}.
We deal with the terms $\omega_1^-$ and $\omega_1^+$ separately (see
equation \eqref{eq-omega1} for the definition of these differentials).

For each edge $e\in E$ from $v$ to $v'$, let 
$\boE_e$ be the space of bounded holomorphic functions
on $\Omega_{v,\epsilon}$ with the sup norm.
This is a Banach algebra for the pointwise product.
If $|p_e^- - \underline{p}_e^-|<\frac{\epsilon}{2}$
then for any $z\in \Omega_{v,\epsilon}$, we have $|z-p_e^{-}|>\frac{\epsilon}{2}$.
Let $f_e(p_e^-)$ be the holomorphic function $\frac{1}{z-p_e^-}$ on
$\Omega_{v,\epsilon}$. The map $p_e^{-}\mapsto f_e(p_e^-)$ is holomorphic and bounded
by $\frac{2}{\epsilon}$ from the disk
$D(\underline{p}_e^-,\frac{\epsilon}{2})$ to $\boE_e$.
Let $\boE$ be the product of the spaces $\boE_e$ for $e\in E$
with the sup norm.
Let $f(\pp)=(f_e(p_e^-))_{e\in E}$.
By lemma \ref{lemma3}, the map $\pp\mapsto f(\pp)$ is smooth from
the ball $B(\underline{\pp},\frac{\epsilon}{2})$ to $\boE$.

We compose this map with
the operator which maps
$(\alpha,f)\in\boB_E\times\boE$
to the differential
$\sum_{e\in E_v^-} \alpha_e f_e dz$ in $\Omega_{v,\epsilon}$, for $v\in V$.
It follows from the definition of the norms that
this operator is bilinear bounded from $\boB_E\times \boE$ to $\boH_1(\Omega_{\epsilon})$.
We deal with $\omega_1^+$ in the same way.
This proves the first point of lemma \ref{lemma-smooth}.

\medskip

The proof of point 2 of lemma \ref{lemma-smooth} is very much similar.
We first deal with the term $\omega_2^-$.
Let $f_{e,n}(p_e^-)\in\boE_e$
be the function $(\frac{\epsilon}{2})^n\frac{1}{(z-p_e^-)^{n}}$. If
$|p_e - \underline{p}_e^-|<\frac{\epsilon}{2}$,
then $||f_{e,n}||_{\infty}\leq 1$.
Define $f(\pp)=(f_{e,n}(p_e^-))_{e\in E,n\in\N}$.
By lemma \ref{lemma3}, the map $\pp\mapsto f(\pp)$ is smooth from
the ball $B(\underline{\pp},\frac{\epsilon}{2})$ to some space $\boE'$,
namely the product of the spaces $\boE_e$ for $e\in E$ and $n\in\N$,
with the sup norm.

We compose this map with the operator
which maps $(\lambda,f)$ to
the differential
$\sum_{e\in E_v^-}\sum_{n\geq 2}
\lambda_{e,n}^- 2^{-n} f_{e,n}^-dz$.
It is straightforward to check that
this operator is bilinear bounded from
$\boB_L\times\boE'$ to $\boH_1(\Omega_{\epsilon})$
(the series converges thanks to the term $2^{-n}$).
We deal with $\omega_2^+$ in the same way.
This proves point 2 of lemma \ref{lemma-smooth}.

\medskip

To prove point 3 of the lemma, we replace the integration circle in the definition of $F_{e,n}^+$
by the fixed circle $C(\underline{p}_e^-,r)$
with $r=\frac{\rho\cc}{2}$.
For each edge $e\in E$,
let $A_e^-$ be the fixed annulus
$\frac{r}{2}<|z-\underline{p}_e^-|<\frac{3r}{2}$.
Let $\boE_e''$ be the spaces of bounded
holomorphic functions on $A_e^-$.
Assume that $||\pp- \underline{\pp}||_{\infty}<\frac{r}{4}$.
If $z\in A_e^-$, then $\frac{r}{8}<|z-p_e^-|<2r$.
Assume that $|t_e|<\frac{\cc\rho^2}{8\ccc}$.
Then for $z\in A_e^-$, we have
$\frac{\ccc|t_e|}{\rho}<|z-p_e^-|<\cc\rho$.
Using the definition of admissible coordinates, we
get $\frac{|t_e|}{\rho}<|z_e^-(z)|<\rho$, so $\varphi_e^+(z)$ is defined.
Let $f_{e,n}\in \boE_e''$
be the function $(\frac{4}{\epsilon})^n(\varphi_e^+(z) - p_e^+)^n$
on $A_e^-$.
By equation \eqref{eq-varphi}, we have
$ |\varphi_e^+ - p_e^+|
\leq \frac{8\ccc^2|t_e|}{\rho\cc}$.
If $|t_e|\leq \frac{\rho\cc\epsilon}{32\ccc^2}$,
then $f_{e,n}$ is bounded by
$1$ on $A_e^-$.
Let $\boE''$ be the product of the spaces $\boE_e''$ for
$e\in E$ and $n\in\N$, with the sup norm.
Let $f=(f_{e,n})_{e\in E,n\in\N}\in\boE''$.
By lemma \ref{lemma3}, the map $(\pp,\tt,\xi)\mapsto f$ is
smooth from
a neighborhood of $(\underline{\pp},0,\underline{\xi})$ to
$\boE''$.

We compose this map with the operator
$$(f,w)\mapsto \disp\left( \int_{C(\underline{p}_e^-,r)}
f_{e,n-1} \,w\right)_{e\in E,n\geq 2}.$$
This operator is bilinear bounded from
$\boE''\times\boH_1(\Omega_{\epsilon})$ to $\boB_L$.
This follows readily from the definition of the norms.
This proves that $F^+$ depends smoothly on parameters. We
deal with $F^-$ in the same way.
This concludes the proof of lemma \ref{lemma-smooth}.
\cqfd
\begin{remark}
\label{remark5}\em
We have used the following fact : if $f(z,w)$ is a bounded holomorphic
function of $z\in\Omega\subset\C$ and $w\in B(0,R)\subset\C^d$,
then the map $w\mapsto f(\cdot,w)$ is holomorphic from the open ball
$B(0,R)$ to the space of holomorphic functions on $\Omega$ with the
sup norm. As in the proof of lemma \ref{lemma3}, this follows from the
fact that we have uniform estimates for the derivatives of $f$ with
respect to $w$.
\end{remark}
\begin{remark}\em
The partial differential of $\omega$ with respect to
$\tt$ at $\tt=0$ can be computed as follows :
let us fix all variables but $\tt$ and use the notation
$\frac{\partial}{\partial\tt}$ for the differential with respect to
$\tt$ at $\tt=0$. We have
$$\omega=\omega_1+\omega_2(F(\tt,\omega)).$$
Take the differential with respect to $\tt$ :
$$\frac{\partial\omega}{\partial \tt}\cdot\hh=\omega_2(\frac{\partial F}{\partial \tt}
(0,\omega)\cdot\hh)+\omega_2(F(0,\frac{\partial\omega}{\partial\tt}\cdot\hh)).$$
The second term is zero because $F=0$ when $\tt=0$.
To compute the first term we write
$$\frac{\partial\varphi_e^+(z)}{\partial t_e}
=\frac{1}{(z_e^+)'(p_e^+)z_e^-(z)}.$$
From this we obtain
$$\frac{\partial\omega}{\partial\tt}\cdot \hh
=\sum_{e\in E_v^+}h_e a_e^+ \frac{dz}{(z-p_e^+)^2}
+\sum_{e\in E_v^-}h_e a_e^- \frac{dz}{(z-p_e^-)^2}
$$
with
$$a_e^+=\frac{-1}{2\pi i (z_e^+)'(p_e^+)}\int_{\partial D_e^-}\frac{\omega}{z_e^-}$$
and a similar definition for $a_e^-$.
In particular, when the coordinates are given as in example \ref{example1} by
$z_e^{\pm}=z-p_e^{\pm}$, we have
$$a_e^+=\frac{-1}{2\pi i}\int_{\partial D_e^-}\frac{\omega}{z-p_e^-}.$$
\end{remark}
\section{Estimate of $\omega$ on the necks}
When $\tt=0$, $\omega_2=0$ so $\omega=\omega_1$ is explicitely
given by formula \ref{eq-omega1}. Theorem \ref{theorem-smooth}
tells us that the restriction of $\omega$ to the domain
$\Omega_{\epsilon}$ depends smoothly on the parameter
$\tt$, so gives us good
control of $\omega$ on $\Omega_{\epsilon}$ for small values of $\tt$
but says nothing outside of $\Omega_{\epsilon}$, i.e. on the necks.

\medskip

In a neighborhood of $p_e^+$, we expect $\omega\simeq
\alpha_e\frac{dz}
{z-p_e^+}$. The following proposition gives a precise
statement.
\begin{proposition}
\label{proposition-necks}
Assume that the degree of vertices is bounded.
There exists a uniform constant $C$ such that for $\tt$
small enough and for any edge $e\in E$
one has
$$\left|\frac{\omega}{dz}-\frac{\alpha_e}{z-p_e^+}\right|
\leq C||\alpha||_{\infty}\qquad
\mbox{ in the annulus $|t_e|^{1/2}\leq |z_e^+|\leq \rho$},$$
$$\left|\frac{\omega}{dz}+\frac{\alpha_e}{z-p_e^-}\right|
\leq C||\alpha||_{\infty}\qquad
\mbox{ in the annulus $|t_e|^{1/2}\leq |z_e^-|\leq \rho$}.$$
\end{proposition}
Proof :
we estimate $\omega_1$ and $\omega_2$ separately.
Consider $z$ in the annulus $|t_e|^{1/2}\leq |z_e^+(z)|\leq\rho$.
Then for $e'\in E_v^+$, $e'\neq e$, one has
$|z-p_{e'}^+|\geq \rho\cc$ and for
$e'\in E_v^-$, one has
$|z-p_{e'}^-|\geq \rho\cc$.
This readily gives the estimate
$$\left|\frac{\omega_1}{dz}-\frac{\alpha_e}{z-p_e^+}\right|
\leq C||\alpha||_{\infty}.$$
By remark \ref{remark2} with $\epsilon=\rho\cc$, we have
$||\omega||_{\infty}\leq C||\alpha||_{\infty}$.
Using that $\lambda=F(\omega)$ and estimate \eqref{eq-estimee-F}, we have
$$\left(\frac{\epsilon}{4}\right)^n|\lambda_{e,n}^+|
\leq Cr \left(\frac{\ccc^2 |t_e|}{r}\right)^{n-1} ||\alpha||_{\infty}.$$
Using the definition of admissible coordinates, we have
$|z-p_e^+|\geq \cc|t_e|^{1/2}$. For $e'\in E_v^+$,
$e'\neq e$, we have 
$|z-p_{e'}^+|\geq \rho\cc\geq\cc|t_e|^{1/2}$. This gives
\begin{eqnarray*}
|\frac{\omega_2^+}{dz}|&\leq& \sum_{e'\in E_v^+}\sum_{n=2}^{\infty}
Cr||\alpha||_{\infty}\left(\frac{\ccc^2|t_e|}{r}\right)^{n-1}\frac{1}
{(\cc |t_e|^{1/2})^n}\\
&\leq& \sum_{e'\in E_v}C||\alpha||_{\infty}\frac{2\ccc^2}{\cc^2}
\end{eqnarray*}
provided $\ccc^2|t_e|^{1/2}\leq \frac{r\cc}{2}$.
We estimate $\omega_2^-$ in the same way. This proves
the first statement of the proposition. The second statement
is similar.
\cqfd
\begin{corollary}
For any $\epsilon\in(0,\cc\rho)$ and any edge $e\in E$, if 
$|\alpha_e|>C\epsilon||\alpha||_{\infty}$,
then $\omega$ has no zero in the annulus bounded in
$\Sigma_{\tt}$ by the circles $C(p_e^-,\epsilon)$
and $C(p_e^+,\epsilon)$, where $C$ is the constant
that appears in proposition \ref{proposition-necks}.
\end{corollary}
Proof : if $z$ is in this annulus, then either
$|z_e^+(z)|\geq |t_e|^{1/2}$ or $|z_e^-(z)|\geq |t_e|^{1/2}$.
Proposition \ref{proposition-necks} gives
$|\frac{\omega}{dz}|>0$.
\cqfd

\section{Decay to a bounded differential}
Given $\alpha\in\ell^{\infty}(E)$, theorem \ref{theorem1} gives
us a differential $\omega\in\boH^1_{\infty}(\Sigma_{\tt})$.
Consider an edge $e\in E$, and assume that we perturb the
parameter $t_e$ while keeping the other parameters
$t_{e'}$ for $e'\neq e$ fixed. This perturbs the differential $\omega$,
and we would like to understand what is the decay of this
perturbation as $v\to\infty$.

\medskip

As in example \ref{example4}, we do this by using weighted
spaces
$\ell^{\infty,\sigma}$
for the parameter $\tt$.
So in this section, we assume that the parameter
$\tt$ is in a neighborhood of some central value $\underline{\tt}$
for the norm $\boB_E$.
There are two observations to be made: firstly, the
central value $\underline{\tt}$ is in $\ell^{\infty}$ but
might not be in $\boB_E$, so we write $\dot{\tt}=
\tt-\underline{\tt}$ and require that $\dot{\tt}\in\boB_E$.
Second, to open node, we need that $||\tt||_{\infty}$ remains small,
so we require that $||\dot{\tt}||_{\infty}\leq C||\dot{\tt}||_{E}$, or
in other words, that $\boB_E\hookrightarrow \ell^{\infty}$.

\medskip

Another motivation is that we will study the smooth dependance
of $\omega$ on $\tt$ for the norm $\boB_E$ instead of the
norm $\ell^{\infty}$ as in section \ref{section-smooth}. This is
more natural for certain applications we have in mind.

\medskip

The same question can be asked for the parameters $\pp$ and
$\xi$.
We define a norm on the space of parameters $\pp$ and $\xi$
(still denoted $||\cdot||_E$ as no confusion can possibly arise) by
$||\pp||_E=||(\max\{|p_e^+|,|p_e^-|\})_{e\in E}||_E$
and $||\xi||_E=||(\max\{||\xi_e^+||,||\xi_e^-||\})_{e\in E}||_E$.
We write $\pp=\underline{\pp}+\dot{\pp}$
and $\xi=\underline{\xi}+\dot{\xi}$
with $\dot{\pp}\in\boB_E$ and $\dot{\xi}\in\boB_E$.
\begin{theorem}
\label{theorem-smooth2}
Assume that $\boB_E\hookrightarrow \ell^{\infty}(E)$.
Then for $||\underline{\tt}||_{\infty}$ small enough,
$\alpha\in\ell^{\infty}(E)$,
and $||\dot{\pp}||_E$, $||\dot{\tt}||_E$ and $||\dot{\xi}||_E$
small enough,
the restriction of
$\omega(\underline{\pp}+\dot{\pp},\underline{\tt}+\dot{\tt},\underline{\xi} +\dot{\xi},\alpha)
-\omega(\underline{\pp},\underline{\tt},\underline{\xi},\alpha)$
to $\Omega_{\epsilon}$ belongs to $\boH^1(\Omega_{\epsilon})$.
Moreover,
the map $(\dot{\pp},\dot{\tt},\dot{\xi},\alpha)\mapsto
\omega(\underline{\pp}+\dot{\pp},\underline{\tt}+\dot{\tt},\underline{\xi} +\dot{\xi},\alpha)
-\omega(\underline{\pp},\underline{\tt},\underline{\xi},\alpha)$
is smooth.
\end{theorem}
The norm on the parameters $\dot{\pp}$, $\dot{\tt}$, $\dot{\xi}$ is
the norm $||\cdot||_E$. The norm on $\alpha$ is the $\ell^{\infty}$ norm.
The norm on the target space is the norm $||\cdot||_{\Omega}$.
Observe that a priori $\alpha\not\in\boB_E$ so
$\omega(\pp,\tt,\xi,\alpha)\not\in\boH^1(\Sigma_{\tt})$, which is
why we have to substract $\omega(\underline{\pp},\underline{\tt},
\underline{\xi},\alpha)$.

\medskip

Returning to the question raised at the beginning of this section,
we take $\boB_E$ to be the space
$\ell^{\infty,\sigma}(E)$ defined in
section \ref{section-decay}.
Fix some edge $e\in E$ and assume that
$\dot{t}_{e'}=0$ for $e'\neq e$.
Then $\dot{\tt}\in\boB_E$, so
$||\omega(\underline{\tt}+\dot{\tt},\alpha)
-\omega(\underline{\tt},\alpha)||_{\infty,\sigma}$ is
finite. In other words, the influence of the parameter $t_e$ on 
the restriction of $\omega$ to $\Omega_{v,\epsilon}$ decays exponentially as
$v\to\infty$.
The same statement holds for the parameters $p_e^{\pm}$ and $\xi_e^{\pm}$.

\medskip

Proof of theorem \ref{theorem-smooth2}. The proof is very similar to
the proof of theorem \ref{theorem-smooth}. We need the following
analogue of lemma \ref{lemma-smooth}.
\begin{lemma}
\label{lemma-smooth2}
For $||\underline{\tt}||_{\infty}$ small enough, $\alpha\in\ell^{\infty}$,
$\lambda\in\ell^{\infty}$, $w\in\boH^1_{\infty}(\Omega_{\epsilon})$, and
$||\dot{\pp}||_{E}$, $||\dot{\tt}||_E$ and $||\dot{\xi}||_E$ small
enough,
\begin{enumerate}
\item $\omega_1(\underline{\pp}+\dot{\pp},\alpha) - \omega_1(\underline{\pp},\alpha)
\in\boH^1(\Omega_{\epsilon})$ and depends smoothly on $\dot{\pp}$,
\item $\omega_2(\underline{\pp}+\dot{\pp},\lambda)- \omega_2(\underline{\pp},\lambda)
\in\boH^1(\Omega_{\epsilon})$ and depends smoothly on $\dot{\pp}$,
\item $F(\underline{\pp}+\dot{\pp},\underline{\tt}+\dot{\tt},\underline{\xi}+\dot{\xi},w)-F(\underline{\pp},\underline{\tt},\underline{\xi},
w)\in\boB_L$ and depends smoothly on $(\dot{\pp},\dot{\tt},\dot{\xi})$.
\end{enumerate}
\end{lemma}
Proof of point 1. Let $\pp=\underline{\pp}+\dot{\pp}$.
We have in $S_v$
$$\omega_1^-(\pp,\alpha)-\omega_1^-(\underline{\pp},\alpha)
=\sum_{e\in E_v^-}\alpha_e\frac{p_e^- - \underline{p}_e^-}{
(z-p_e^-)(z-\underline{p}_e^-)}dz
=\sum_{e\in E_v^-}\alpha_e \dot{p}_e^- g_e dz$$
with $g_e=\frac{1}{(z-p_e^-)(z-\underline{p}_e^-)}$.
Let $g=(g_e)_{e\in E}$. As in the proof of point 1 of lemma \ref{lemma-smooth},
the map $\pp\mapsto g$ is smooth from the ball
$||\pp-\underline{\pp}||_{\infty}<\frac{\epsilon}{2}$
to $\boE$. Since $\boB_E\hookrightarrow \ell^{\infty}$, the map
$\dot{\pp}\mapsto g$ is smooth in a neighborhood of $0$. We compose this map
with the operator
$$(\alpha,g,\dot{\pp})\mapsto \sum_{e\in E_v^-}\alpha_e g_e \dot{p}_e^- dz.$$
This operator is trilinear bounded from
$\ell^{\infty}(E)\times\boE\times\boB_E$ to $\boH^1(\Omega_{\epsilon})$.
We deal with $\omega_1^+$ in the same way.

\medskip

Proof of point 2. We write in $S_v$
$$\omega_2^-(\pp,\lambda)-\omega_2^-(\underline{\pp},\lambda)
=\sum_{e\in E_v^-}\lambda_{e,n}^- 2^{-n} g_{e,n} \dot{p}_e^- dz$$
with
$$g_{e,n}=\frac{f_{e,n}(p_e^-) - f_{e,n}(\underline{p}_e^-)}
{p_e^- - \underline{p}_e^-}$$
and $f_{e,n}$ as in the proof of point 2 of lemma \ref{lemma-smooth}.
Now if a holomorphic function $f(z)$ is bounded by $C$ on the disk $D(0,R)$,
then its derivative is bounded by $\frac{2C}{R}$ on the disk $D(0,\frac{R}{2})$
by Cauchy estimate. Hence by the mean value inequality,
the function $\frac{f(z)-f(w)}{z-w}$ is bounded by $\frac{2C}{R}$ on
$D(0,\frac{R}{2})\times D(0,\frac{R}{2})$.
Hence, from the bound $|f_{e,n}|\leq 1$ for
$|p_e^- - \underline{p}_e^-|<\frac{\epsilon}{2}$,
we get a uniform bound of $g_{e,n}$ for
$|p_e^- - \underline{p}_e^-|<\frac{\epsilon}{4}$.
Then we conclude as in the proof of point 1 by composing with a bounded
trilinear operator.

Proof of point 3. We first prove that
$F(\pp,\tt,\xi,w)-F(\underline{\pp},\tt,\xi,w)$ depends smoothly on
$(\dot{\pp},\tt,\xi)$, hence on $(\dot{\pp},\dot{\tt},\dot{\xi})$
since $\boB_E\hookrightarrow \ell^{\infty}$.
To do this, we adapt the proof of point 3 of lemma \ref{lemma-smooth}
exactly as we did above.
Then we prove in the same way that $F(\underline{\pp},\tt,\xi,w)-F(\underline{\pp},\underline{\tt},
\xi,w)$ depends smoothly on $(\dot{\tt},\xi)$, and
finally that $F(\underline{\pp},\underline{\tt},\xi,w)-F(\underline{\pp},\underline{\tt},\underline{\xi},w)$ depends smoothly on $\dot{\xi}$.
\cqfd

\medskip

Proof of theorem \ref{theorem-smooth2}.
Let
$$G(\pp,\tt,\xi,\alpha,\lambda)=F(\pp,\tt,\xi,\omega_1(\pp,\alpha)+
\omega_2(\pp,\lambda)).$$
Recall that $G$ is contracting with respect to $\lambda$ and
that we found $\lambda$ as a fixed point
of $\lambda\mapsto G(\pp,\tt,\xi,\alpha,\lambda)$.
Let $\underline{\lambda}$ be the fixed point of 
$\lambda\mapsto G(\underline{\pp},\underline{\tt},\underline{\xi},\alpha,\lambda)$,
so that
$\omega(\underline{\pp},\underline{\tt},\underline{\xi},\alpha)=
\omega_1(\underline{\pp},\alpha)+\omega_2(\underline{\pp},\underline{\lambda})$.
Write $\lambda=\underline{\lambda}+\dot{\lambda}$ and define
$$H(\dot{\pp},\dot{\tt},\dot{\xi},\alpha,\dot{\lambda})
=G(\pp,\tt,\xi,\alpha,\lambda)-G(\underline{\pp},\underline{\tt},\underline{\xi},
\alpha,\underline{\lambda}).$$
Then $\dot{\lambda}$ is a fixed point of $\dot{\lambda}\mapsto
H(\dot{\pp},\dot{\tt},\dot{\xi},\alpha,\dot{\lambda})$.
By claim \ref{claim-smooth2} below and the fixed point theorem with
parameters, $\dot{\lambda}\in\boB_L$ and depends smoothly on
$(\dot{\pp},\dot{\tt},\dot{\xi})$.
Then we write
$$\omega(\pp,\tt,\xi,\alpha)-\omega(\underline{\pp},\underline{\tt},
\underline{\xi},\alpha)
=(\omega_1(\pp,\alpha)-\omega_1(\underline{\pp},\alpha))
+(\omega_2(\pp,\lambda)-\omega_2(\underline{\pp},\lambda))
+\omega_2(\underline{\pp},\dot{\lambda}).$$
It follows from lemmae \ref{lemma-smooth} and \ref{lemma-smooth2} that
all three terms are in $\boH^1(\Omega_{\epsilon})$ and depend smoothly on
$(\dot{\pp},\dot{\tt},\dot{\xi})$.
\begin{claim}
\label{claim-smooth2}
For $||\underline{\tt}||_{\infty}$ small enough,
$(\dot{\pp},\dot{\tt},\dot{\xi})$ in a neighborhood of $0$
and $\dot{\lambda}\in\boB_L$, we have $H(\dot{\pp},\dot{\tt},\dot{\xi},\alpha,
\dot{\lambda})\in \boB_L$ and it depends smoothly on
$(\dot{\pp},\dot{\tt},\dot{\xi})$.
Moreover, it is contracting with respect to $\dot{\lambda}$.
\end{claim}
Proof. We write
\begin{eqnarray*}
H(\dot{\pp},\dot{\tt},\dot{\xi},\alpha,\dot{\lambda})
&=& F(\pp,\tt,\xi,\omega_1(\pp,\alpha))-F(\underline{\pp},\underline{\tt},\underline{\xi},\omega_1(\underline{\pp},\alpha))\\
& & + F(\pp,\tt,\xi,\omega_2(\pp,\lambda))-F(\underline{\pp},\underline{\tt},
\underline{\xi},\omega_2(\underline{\pp},\underline{\lambda})).
\end{eqnarray*}
The first conclusion then follows from lemma \ref{lemma-smooth2}.
To prove $H$ is contracting, we write
$$H(\dot{\pp},\dot{\tt},\dot{\xi},\alpha,\dot{\lambda})
-H(\dot{\pp},\dot{\tt},\dot{\xi},\alpha,0)
=F(\pp,\tt,\xi,\omega_2(\pp,\dot{\lambda})).$$
We have already seen in section \ref{section-proof}
that this operator is contracting with respect to 
$\dot{\lambda}$.
This proves the claim and theorem \ref{theorem-smooth2}.
\cqfd
\section{Meromorphic differentials}
\label{section-meromorphic}
In this section we adapt the results of the previous sections to the case
of meromorphic 1 forms.
\medskip

On a compact Riemann surface, one can define a meromorphic 1-form $\omega$ by
prescribing its poles and principal parts at each pole, with the only
condition that the sum of the residues be zero. Morevoer, this defines
$\omega$ uniquely up to a holomorphic 1-form.

The definition of a meromorphic differential on a Riemann surface with nodes
is the same as that of a regular differential, except that it is allowed to
have poles away from the nodes.

In each sphere $S_v$, we select $n_v$ distincts points $q_{v,i}$,
$1\leq i\leq n_v$, distinct from the nodes,
to be the poles of $\omega$.
One of these points may be $\infty$.
For each pole $q_{v,i}$, we choose an integer $m_{v,i}\geq 1$, to be
the order of the pole.
We define the divisor $Q$ as the formal sum $\disp\sum_{v\in V}
\sum_{i=1}^{n_v} m_{v,i} q_{v,i}$.
We write $(\omega)\geq Q$ to say that $\omega$ has at most poles of order
$m_{v,i}$ at $q_{v,i}$, and is otherwise regular (so it may also have simple poles
at the nodes as a regular differential).

If $\omega$ is such a differential, we may write its principal part at $q_{v,i}$
as follows if $q_{v,i}\neq \infty$ :
$$P_{v,i}(\omega)=\sum_{n=1}^{m_{v,i}} a_{v,i,n} \frac{dz}{(z-q_{v,i})^n}.$$
By definition this means that $\omega - P_{v,i}(\omega)$ is holomorphic in
a neighborhood of $q_{v,i}$.
If $q_{v,i}=\infty$, then we use $w=\frac{1}{z}$ as a coordinate in a neighborhood
of $\infty$ to write the principal part :
$$P_{v,i}(\omega)=\sum_{n=1}^{m_{v,i}} a_{v,i,n} \frac{dw}{w^n}
=-\sum_{n=1}^{m_{v,i}} a_{v,i,n} z^{n-2} dz.$$
The residue of $\omega$ at $q_{v,i}$ is $a_{v,i,1}$.
If we let $\alpha_e=\int_{\gamma_e}\omega$, we have by the residue theorem :
\begin{equation}
\label{eq-residue}
\forall v\in V,\qquad \sum_{e\in E_v^+}\alpha_e-\sum_{e\in E_v^-}\alpha_e
+2\pi i\sum_{i=1}^{n_v} a_{v,i,1}=0.
\end{equation}
Essentially, this is the only obstruction to define a meromorphic differential by prescribing its
principal parts and periods.

We write $P=(P_{v,i})_{v\in V,1\leq i\leq n_v}$.
We define $\boB_P$ to be the space of principal parts $P$ such that the following
norm is defined :
$$||P||_{P}=\left\|(\sup
\{|a_{v,i,n}| \; :\;
1\leq i \leq n_v,\,1\leq n\leq m_{v,i}\}
)_{v\in V}\right\|_V.$$
\begin{hypothesis}
\label{hypothesis-meromorphic} We assume that
\begin{enumerate}
\item the coordinates and the norms are admissible,
\item the number of poles in each sphere, and the orders of the poles, are uniformly bounded,
\item there exists $r>0$ such that for each vertex $v\in V$,
the disks $D_e^-$ for $e\in E_v^-$, $D_e^+$
for $e\in E_v^+$ and $D(q_{v,i},r)$ for $1\leq i\leq n_v$ are disjoint.
\end{enumerate}
\end{hypothesis}
In case $q_{v,i}=\infty$, the disk $D(q_{v,i},r)$ should be understood as
the disk $|w|\leq  r$, or equivalently, $|z|\geq \frac{1}{r}$.
In this setup, the definition of $\Omega_{v,\epsilon}$ must be changed as follows :
Given some positive $\epsilon$ less than $r$, we define 
the domain $\Omega_{v,\epsilon}$ as the Riemann sphere $S_v$ minus the disks
$D(p_e^-,\epsilon)$
for $e\in E_v^-$, $D(p_e^+,\epsilon)$ for $e\in E_v^+$ and
$D(q_{v,i},\epsilon)$ for $1\leq i\leq n_v$.
We define $\Omega_{\epsilon}$
as the disjoint union of the domains $\Omega_{v,\epsilon}$ for $v\in V$.
The norm $||\cdot||_{\Omega}$ is then defined as in section \ref{section-admissiblenorms}
with this new definition of $\Omega_{v,\epsilon}$.

Let $\boM^1(\Sigma_{\tt},Q)$ be the space of meromorphic differentials $\omega$ on $\Sigma_{\tt}$
such that $(\omega)\geq Q$ and the norm $||\omega||_{\Omega}$ is defined.
The following theorem is the generalisation of theorem \ref{theorem2} to the case of meromorphic
differentials.
\begin{theorem}
\label{theorem-meromorphic}
For $\tt$ small enough,
the operator $\omega\mapsto ((\int_{\gamma_e}\omega)_{e\in E},P(\omega))$ is an isomorphism from
$\boM^1(\Sigma_{\tt},Q)$ to the subspace of vectors in $\boB_E\times\boB_P$ which satisfy the
condition \eqref{eq-residue}.
\end{theorem}
Proof :
the proof of this theorem is exactly the same as the proof of theorem \ref{theorem2},
except that in point 2 of the proof,
the definition of $\omega_1$, equation \eqref{eq-omega1},
must be changed to take into account the poles.
Given $\alpha\in\boB_E$ and $P\in\boB_P$, satisfying \eqref{eq-residue},
we are asked to construct a meromorphic differential
$\omega$ with prescribed periods $\alpha$ and prescribed
principal parts $P$.
If the point $\infty$ in $S_v$ is not a pole, we define $\omega_1$ in
$S_v$ as
$$\omega_1=
\frac{1}{2\pi i}\sum_{e\in E_v^+} \alpha_e\frac{dz}{z-p_e^+} 
-\frac{1}{2\pi i}\sum_{e\in E_v^-} \alpha_e\frac{dz}{z-p_e^-} 
+\sum_{i=1}^{n_v}\sum_{n=1}^{m_{v,i}} a_{v,i,n} \frac{dz}{(z-q_{v,i})^n}.$$
Equation \eqref{eq-residue} guarantees that $\omega_1$ is holomorphic at $\infty$.
If the point $\infty$ in $S_v$ is a pole, say $q_{v,1}=\infty$,
the definition of $\omega_1$ in $S_v$ is
the same except that the term corresponding to $i=1$ in the last sum is
replaced by
$$-\sum_{n=2}^{m_{v,i}} a_{v,i,n} z^{n-2} dz.$$
Mind the fact that the sum starts at $n=2$ and not $n=1$, so that this term 
is holomorphic at $0$. Equation \eqref{eq-residue} ensures that the residue of
$\omega_1$ at $\infty$ is $a_{v,1,1}$, so $\omega_1$ has the required principal part
at infinity.

With this modification, point 1 of lemma \ref{lemma1} is replaced by
$$||\omega_1||_{\Omega}\leq \frac{C}{\epsilon}||\alpha||_E+\frac{C}{
\epsilon^k}||P||_P$$
where $k=\max\{|m_{v,i}|\,:\, v\in V,1\leq i\leq n_v\}$.
The proof of this fact is straightforward. The rest of the proof is unchanged.
\cqfd
\begin{example} \em Choose a vertex $v_0$. There exists a unique meromorphic differential on $\Sigma_{\tt}$
which has just one double pole at $\infty$ in $S_{v_0}$, with principal part $dz$, and has vanishing periods
on all cycles $\gamma_e$, $e\in E$ (a normalised differential of the second kind).
We can use weighted $\ell^{\infty}$ norms to study its decay as in section \ref{section-decay}, with
the same conclusion.
\end{example}
\begin{example}
\em Choose two vertices $v_1$ and $v_2$. Let $\gamma$ be a path from $v_1$ to $v_2$. Orient $\gamma$ from $v_1$ to
$v_2$ and orient the rest of $\Gamma$ arbitrarily. Choose a point $q_1$ in $S_{v_1}$ and a point $q_2$
in $S_{v_2}$. Then there exists a unique meromorphic differential which has simple poles at $q_1$ and
$q_2$, with respective residues $1$ and $-1$, and has period on the cycles $\gamma_e$ equal to
$2\pi i$ if the edge $e$ belongs to $\gamma$ and zero otherwise (a normalised differential of the third kind).
\end{example}
\begin{remark} All the other results in this paper have a natural generalisation
to the case of meromorphic differentials. For instance, the generalisation
of theorem \ref{theorem-smooth} would be that the restriction of $\omega$
to $\Omega_{\epsilon}$
depends smoothly on the parameters $(\pp,\tt,\xi,\qq)$ in
a neighborhood of $(\underline{\pp},0,\underline{\xi},\underline{\qq})$ in $\ell^{\infty}$
norm. It suffices, in the proofs, to replace $\omega_1$ by the above definition.
The reason we did not consider meromorphic differentials from the very beginning
is that there are more parameters in this case, so altough conceptually the same,
the proofs would have been longer.
\end{remark}

\bigskip

\noindent
Martin Traizet\\
Laboratoire de Math\'ematiques et Physique Th\'eorique\\
Universit\'e de Tours\\
37200 Tours, France\\
\verb$martin.traizet@lmpt.univ-tours.fr$


\begin{thebibliography}{9}
\bibitem{brown-pearcy} A. Brown, C. Pearcy : Introduction to Operator Theory I - Elements of Functional Analysis. Springer Verlag (1977).
\bibitem{fay} J. D. Fay : Theta Functions on Riemann Surfaces. Springer Verlag, Lecture notes in mathematics No 352 (1973).
\bibitem{imayoshi} Y. Imayoshi, M. Taniguchi: An introduction to Teichmuller spaces. {\em Springer Verlag} (1992)
\bibitem{masur} H. Masur : The extension of the Weil Petersson metric to the boundary of Teichmuller space. Duke Math. J. 43 (1976), 623--635.
\bibitem{filippo} F. Morabito, M. Traizet : Non-periodic Riemann
examples with handles. {\em Preprint} (2010)
\bibitem{traizet} M. Traizet : An embedded minimal surface with no symmetries.  J. Differential Geometry 60:1 (2002), 103--153,
\bibitem{soo}Soo Bong Chae : Holomorphy and calculus in normed spaces. 
\end{thebibliography}
\end{document}